\documentclass[preprint,12pt]{elsarticle}

\usepackage[utf8]{inputenc}
\usepackage{eqnarray,amsmath}
\usepackage{amsfonts,amsthm}
\usepackage{graphicx}
\usepackage{hyperref}    
\hypersetup{
    colorlinks,
    linktocpage,
    citecolor=blue,
    filecolor=black,
    linkcolor=blue,
    urlcolor=magenta,
}
\usepackage{wrapfig}
\usepackage{multirow}
\usepackage{caption}
\usepackage{booktabs}
\usepackage{multicol}

\usepackage{amssymb}
\usepackage{mathrsfs}
\usepackage{mathtools}
\usepackage{xspace}
\usepackage[stable]{footmisc}
\usepackage{wrapfig}
\usepackage{stackengine,xcolor}
\usepackage{booktabs} 
\usepackage{comment}
\usepackage{epstopdf}
\usepackage{float}
\usepackage[font=normal]{subcaption}
\usepackage{stmaryrd}
\usepackage{subfiles}
\usepackage[margin=2.5cm]{geometry}
\usepackage[T1]{fontenc}
\usepackage{babel}
\usepackage{pifont}

\usepackage{bigints}
\usepackage{placeins}
\usepackage[wby]{callouts}
\usepackage{tikz}
\usepackage{subcaption}
\usetikzlibrary {arrows.meta,bending} 
\usepackage[ruled,vlined]{algorithm2e}
\SetKwInput{KwInput}{Input}                
\SetKwInput{KwOutput}{Output}

\usepackage{lineno}

\newcommand{\xx}{\mathbf{x}}

\newcommand{\rmd}{\mathrm{d}}

\usepackage{cleveref}
\usepackage{blindtext}

\graphicspath{{./Figures/}}

\newcommand{\answer}[1]{{\color{black}{#1}}}
\begin{document}

\begin{frontmatter}




\title{\answer{A Generalized Finite Difference-Based Fragile Points Method for Heat Conduction Problems in Non-Homogeneous Media}}

\author[iitm]{Smriti}
\author[iitm]{Sundararajan Natarajan\corref{cor}\fnref{fnlab1}\tnoteref{label2}}
\address[iitm]{Department of Mechanical Engineering, Indian Institute of Technology Madras, Chennai 600036, India.}
\author[iitm]{Abinash Malla}

\tnotetext[label2]{Pandurangan Faculty Fellow}
\cortext[cor]{Corresponding author}
\fntext[fnlab1]{Department of Mechanical Engineering, Indian Institute of Technology Madras, Chennai - 600036, Tamil Nadu, India. Email: snatarajan@iitm.ac.in}


\begin{abstract}
This paper presents an enhanced formulation of the Fragile Points Method (FPM), a truly meshless approach for efficiently modeling implicit interfaces in two-dimensional differential equations involving non-homogeneous materials. The proposed framework eliminates the need for specialized numerical integration techniques and provides a systematic mathematical foundation for solving interface problems. Discontinuities in both primary and secondary variables across interfaces are naturally handled through the inherently discontinuous shape functions of FPM. Unlike conventional Galerkin methods, FPM employs simple, local, point-based polynomial trial and test functions constructed via a generalized finite difference approach. These discontinuous functions bypass the continuity requirements of standard Galerkin frameworks. To address the resulting inconsistency due to discontinuities, we incorporate numerical flux corrections inspired by the discontinuous Galerkin method. The proposed method is validated through several benchmark problems, demonstrating its efficiency and robustness.
\end{abstract}


\end{frontmatter}

\section{Introduction}
\answer{Interfaces are ubiquitous in nature and serve as a critical boundary that facilitates transfer of forces between components. The physical processes occurring at the interface is typically modeled using non-homogeneous second order elliptic equations.
These models are essential for providing an accurate description of systems in which discrete properties across an interface play a major role in determining the overall behaviour. Common applications include heat transfer across various media, diffusion processes in heterogeneous materials, and groundwater flow through distinct soil layers. Elliptic interface problems pose challenges due to discontinuities at the interfaces between different materials, leading to abrupt changes in properties like conductivity or permeability. Accurately capturing these changes is essential for reliable simulation results. Solving these problems is critical across numerous fields such as biomedical engineering, fluid dynamics (where different fluids interact), electronic device optimization (where different conducting materials affect the flow of electric current), the design of composite materials with desired properties, environmental science, and more. }

Typically numerical methods are used for tackling these problems, as analytical solutions are
often unattainable. Commonly used numerical methods include the finite element method (FEM)\cite{hansbo2002unfitted}, the scaled boundary finite element method (SBFEM)\cite{dsouza2021robust}, and the extended finite element method (XFEM) combined with the immersed interface method (IIM)\cite{vaughan2007comparison}. These techniques discretize the domain into smaller nonoverlapping segments, solving the governing equations while managing discontinuities at the interfaces. Adaptive mesh refinement is frequently utilized near interfaces to capture these discontinuities with greater accuracy. Recently, Atluri et al. \cite{paper1} introduced a truly meshfree method that employs point-wise simple polynomial trial and test functions through generalized finite difference methods, known as the fragile points method (FPM). Since these functions are defined point-wise, both the test and trial functions are discontinuous at the cell boundaries. To enhance consistency, and inspired by the discontinuous Galerkin approach, numerical flux corrections are employed \cite{paper2}. The salient features of this framework include: (a) trial and test functions that are local to the cell and inherently discontinuous, making the framework ideal for modeling problems with strong and weak discontinuities; (b) the use of simple polynomial functions, which eliminates the need for specialized numerical integration techniques; and (c) the ability to efficiently manage interface problems and interface boundary conditions without a significant computational burden.
\newline

The accuracy, consistency in the solution, and robustness of the FPM when applied to 1D and 2D Poisson's equations were studied in \cite{paper1}. Additionally, the FPM has been utilized to investigate transient heat conduction problems in \cite{paper3,paper5}, employing the generalized finite difference technique to derive trial and test functions pointwise. Both the primal FPM and mixed FPM approaches have been used to analyze flexoelectric effects in 2D dielectric materials \cite{paper7, Guan2020ANM}. Yang et al. \cite{paper4} developed the FPM for linear elasticity, performing convergence studies on benchmark elasticity problems and simulating the crack propagation and initiation problems. In the context of cardiac electrophysiology, Mountris et al. \cite{mountris2021meshfree} used the FPM to derive the cardiac monodomain model  and estimated ventricular and atrial fiber configurations for four distinct geometries \cite{mountris2022meshless}. Furthermore, Mountris et al. \cite{mountris2023explicit} formulated an explicit Lagrangian algorithm using the FPM to simulate the deformation of hyperelastic materials undergoing large deformation. The FPM has also been applied to study interface debonding in U-notched structures \cite{wang2022fragile} and has proven effective in addressing the two-dimensional hyperbolic telegraph equation \cite{paper6}, as well as linear and non-linear wave equations in cracked domains \cite{haghighi2023study}.
\newline

As seen from the literature, since its inception, the FPM has been utilized to address a variety of problems involving both linear and non-linear equations. The method's robustness and accuracy have been validated for both smooth and singular problems. However, to the best of the authors' knowledge, the FPM has not yet been applied to model implicit interfaces in two-dimensional differential equations within non-homogeneous materials, which are crucial in many engineering disciplines. This work makes several important contributions:
\begin{itemize}
    \item[\ding{113}] \textbf{Extension of FPM to Implicit Interface Problems:} This study extends the FPM to model implicit interfaces in non-homogeneous materials, addressing a notable gap in the existing literature.
    \item[\ding{113}] \textbf{Development of a Comprehensive framework:} A complete mathematical framework has been established to handle implicit/explicit interfaces, making it a promising alternative for both the finite element-based approach and the mesh-free methods. 
    \item[\ding{113}] \textbf{Efficient Numerical Integration:} Proposed method avoids the need for specialized integration techniques commonly required in mesh-free methods and allows the construction of interpolation shape functions without much computational complexity. Furthermore, a simple numerical integration technique is employed to compute the domain integrals.
\end{itemize}
 The structure of the remainder of this paper is organized as follows: \Cref{gov eq} details the governing differential equation for heat conduction in an anisotropic non-homogeneous medium.  \Cref{fpm overview} outlines the implementation of the Fragile Points Method (FPM) for this equation. A discussion on numerical flux correction for ensuring stability and consistency in the solutions, is presented in \Cref{numflux}.  The accuracy, stability and robustness of the FPM in addressing the global response of different interface problems are demonstrated in \Cref{numexap}, and the paper concludes with final remarks in the last section.

\section{Governing equations} \label{gov eq}
\begin{figure}[H]
    \centering
    \includegraphics[scale=0.5]{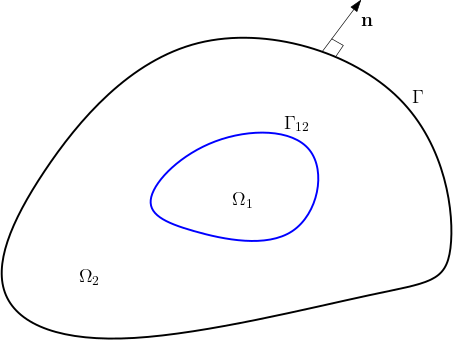}
    \caption{Solid body occupying the domain $\Omega=\Omega_1 \cup \Omega_2$ with inclusion $\Omega_1$, matrix $\Omega_2$ and an interface $\Gamma_{12}$. The external boundary $\Gamma$ comprises Dirichlet $(\Gamma_D)$ and Neumann $(\Gamma_N)$ boundaries, with $\mathbf{n}$ as the outward unit normal to $\Gamma$.}
    \label{fig:inter}
\end{figure}

\noindent In this paper, we investigate heat conduction in an anisotropic heterogeneous medium, using it as a case study to demonstrate the robustness of the Fragile Points Method (FPM) to handle interface problems. Consider a solid body occupying a domain $\Omega \subset \mathbb{R}^2$, enclosed by an external boundary $\Gamma$. This boundary is partitioned into two disjoint segments: the Dirichlet boundary $\Gamma_D \subseteq \Gamma$ and the Neumann boundary $\Gamma_N \subseteq \Gamma$, such that $\Gamma_D \cap \Gamma_N = \emptyset$ and $\Gamma_D \cup \Gamma_N = \Gamma$. The domain $\Omega$ comprises two distinct material regions: an inclusion $\Omega_1$ and a surrounding matrix $\Omega_2$, satisfying $\Omega = \Omega_1 \cup \Omega_2$ and $\Omega_1 \cap \Omega_2 = \emptyset$. These subdomains are separated by an internal interface denoted by $\Gamma_{12}$. See \Cref{fig:inter} for a schematic description. The combined boundary consisting of the external boundary $\Gamma$ and the internal interface $\Gamma_{12}$ is denoted by $\tilde{\Gamma} = \Gamma \cup \Gamma_{12}$.
\noindent The steady state heat conduction within this domain is governed by the following set of equations:
\begin{equation}\label{eqn:thermal_govern}
\left\{
\begin{aligned}
\boldsymbol\nabla \cdot (\beta(\xx) \boldsymbol\nabla u(\xx)) =Q(\xx) \quad {\rm on}  \quad \xx \in \Omega, \\
u(\xx)=g_D(\xx) \quad {\rm in}  \quad \xx \in  \Gamma_D,\\
\beta(\xx) \boldsymbol\nabla u \cdot \mathbf{n}= {q} \quad {\rm in}  \quad \xx \in  \Gamma_N, \\
\llbracket u \rrbracket = 0 \quad {\rm in}  \quad \xx \in  \Gamma_{12},\\
\llbracket \beta \boldsymbol\nabla u \cdot \mathbf{n} \rrbracket =0 \quad {\rm in}  \quad \xx \in  \Gamma_{12},
\end{aligned}
\right.
\end{equation}
where $\boldsymbol\nabla=\{\frac{\partial}{\partial x}\quad \frac{\partial}{\partial y}\}^{\rm T}$ is the gradient operator, $u(\xx)$ is the temperature distribution at a point $\xx$ of the body, $\mathbf{n}$ denotes the outward pointing unit normal to $\Gamma$, the heat flux is denoted by ${q}$ and the internal heat source within the domain is represented by $Q(\xx)$. The diffusion coefficient of the material denoted by $\beta(\xx)$ is a piecewise constant function defined as:
\begin{equation}
   \beta(\xx) = \left\{
\begin{array}{ll} 
\beta_1 \quad {\rm in}  \quad \Omega_1 \\
\beta_2 \quad {\rm in}  \quad \Omega_2
\end{array}
    \right.
\end{equation}
To facilitate numerical solutions, we assume the trial and test functions are selected from function spaces $\mathcal{U}_h$ and $\mathcal{V}_h$ defined as: 
\begin{equation}
\left\{
\begin{aligned}
  & \mathcal{U}_h \subset \mathcal{U}=\{u(\xx) \in H^1(\Omega): u=g_D(\xx) \,\, \text{on}\,\, \Gamma_{D}\}, \\
  &  \mathcal{V}_h \subset \mathcal{V}=\{v(\xx) \in H^1(\Omega): v=0\,\, \text{on}\,\, \Gamma_{D}\}.
\end{aligned} 
\right.
\end{equation}
The domain $\Omega$ is partitioned into non-overlapping element $\Omega_{h}$. The unknown temperature field is approximated by finite-dimensional piecewise function  $u_h(\xx)= \sum_{a} N_{a}(\xx)u_{a}$, where $N_{a}(\xx)$ represents the nodal basis function and $u_a$ is the nodal unknown temperature. Following the standard Bubnov-Galerkin approach, the weak form of the heat conduction equation is expressed as: Find $u \in H^1$ such that $\forall v \in \mathcal{V}_h$,
\begin{equation}\label{eqn: weak_gov EB1}
    \int_{\Omega}\boldsymbol\nabla{v}^{T}\beta\boldsymbol\nabla u ~\rmd\Omega = \int_{\Omega}Q v ~\rmd\Omega + \int_{ \tilde{\Gamma}}v \mathbf{n}^{T}\beta \boldsymbol\nabla{u} ~\text{d}\Gamma.
\end{equation}
The shape functions, denoted as $N_a(\xx)$, are usually continuous within each finite element but may exhibit discontinuities across element boundaries. These functions are typically chosen to be piecewise polynomials within each element of the discretized domain $\Omega_h$ ensuring that the approximation $u_h(\xx)$
inherits the continuity and smoothness properties of the shape functions. The discontinuities across element boundaries are physically meaningful, often corresponding to material interfaces or jumps in properties.

\section{Point-based discontinuous trial and test functions for Fragile Points method}\label{fpm overview}
\begin{figure}[H]
\centering
\includegraphics[width=0.55\linewidth]{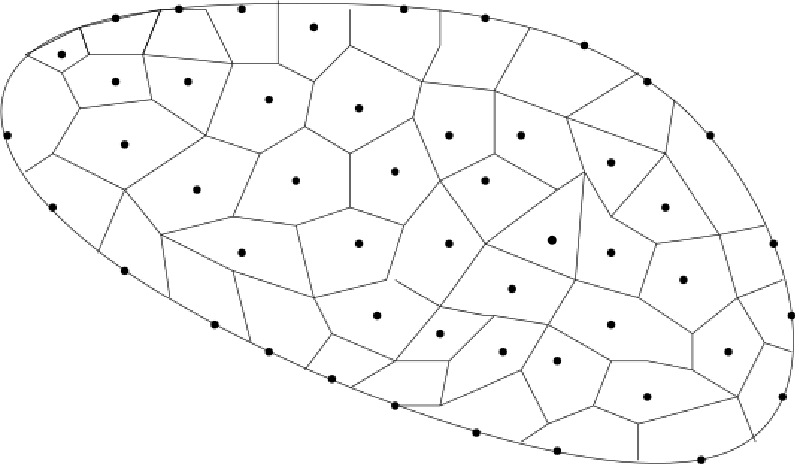}
\caption{Domain and its partition}
\label{fig:figg}
\end{figure}

\begin{figure}[!htb]
    \centering
    \begin{minipage}{.45\textwidth}\hspace{-.7cm}
    \includegraphics[scale=0.5]{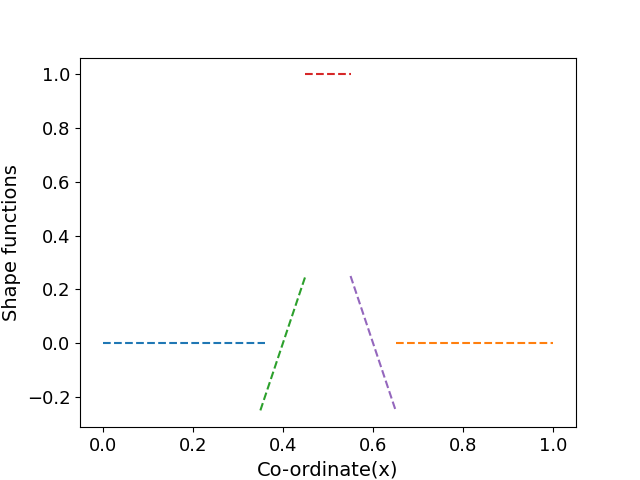}
    \captionsetup{justification=centering}
    \caption{1-D domain}
    \label{fig:shpfunc1}     
    \end{minipage}\hfill
\begin{minipage}{.45\textwidth}\hspace{-0.8cm}\vspace{-.5cm}
    \includegraphics[scale=0.7]{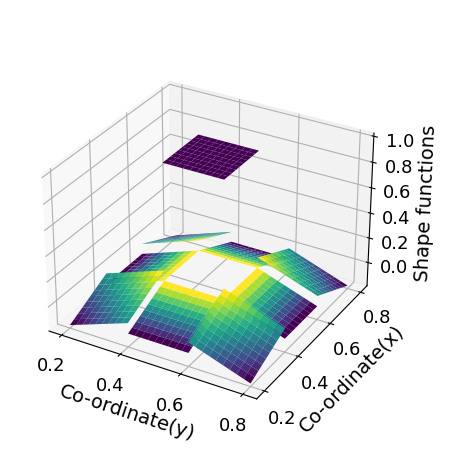}
    \captionsetup{justification=centering}
    \caption{ 2-D domain}
    \label{fig:shpfunc2}
    \end{minipage}
    \end{figure}

\noindent One of the fundamental aspects of the Fragile Points Method (FPM) is its ability to handle problems with complex geometries and potential material discontinuities. This is achieved through the use of local, discontinuous trial and test functions. The first step involves introducing a set of random points within the domain. Then, with these points, the domain can be partitioned into conforming and non-overlapping subdomains in such a way that each subdomain contains only one point. This characteristic grants FPM significant flexibility in handling complex geometries, as the subdomain shape can adapt to the point distribution and can be of any geometry as shown in \Cref{fig:figg}. This advantage arises because the trial and test functions are constructed solely based on the random points within the domain, and their form is entirely independent of the specific subdomain shapes. For creating subdomains, many partitioning schemes can be employed, e.g., the Voronoi Diagram partition, as well as quadrilateral, and triangular partition, even FEM elements, can be converted into FPM subdomains by deﬁning the centroid of each FEM geometrical element as the internal Fragile Point.
\noindent Then, within each subdomain, a simple,
local, discontinuous polynomial trial and test function is established. The trial function, denoted as $u_h(\xx)$, approximates the unknown solution at any point $\mathbf{x}$ located within subdomain $E_0$. This subdomain also houses an internal point, designated as $P_0$. The polynomial trial function can be defined in terms of the values of $u$ and its gradients at $P_0$. For simplicity, we use 2D linear functions. To construct the trial function, we leverage a Taylor series expansion centered at $P_0$. Mathematically, such a trial function is expressed as:
      \begin{align}
          u_{h}(\mathbf{x})=u_{0}+(\mathbf{x}-\mathbf{x_{0}})\cdot\boldsymbol\nabla{u}|_{P_0},\qquad \mathbf{x}\in E_0,
          \label{eqn:trialeq}
      \end{align}
where $u_{0}$ represents the value of the solution field at the internal point $P_0$, $\mathbf{x}$ denotes the position vector of any point within subdomain $E_0$, $\mathbf{x}_0$
  is the position vector of $P_0$ and $\boldsymbol\nabla u|_{P_0}$
  symbolizes the gradient of the solution field, evaluated at $P_0$.
\begin{figure}[!htb]
          \centering
           \subfloat{\includegraphics[width=0.65\textwidth]{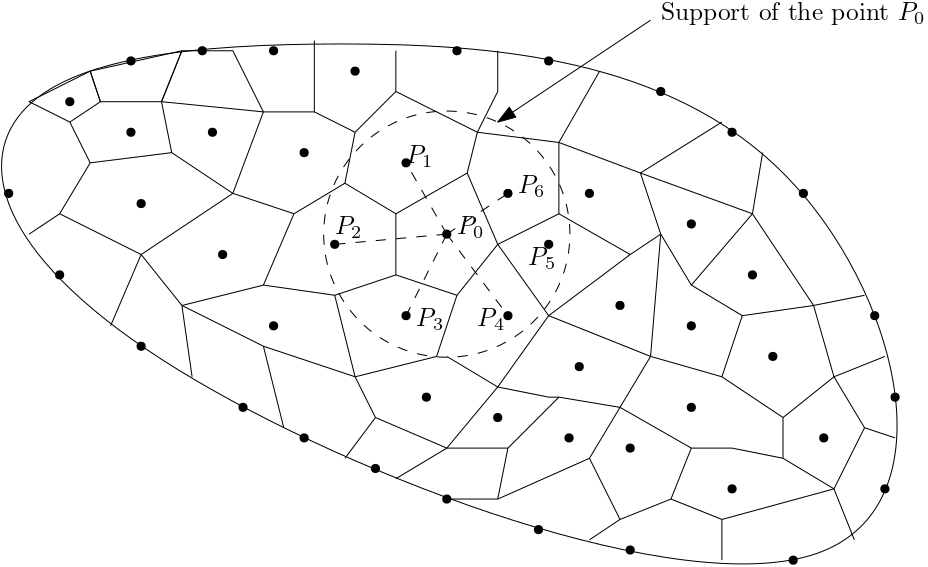}\label{fig:f1}}
           \hfill
            \caption{Domain with support for point $P_{0}$ }
            \label{fig:support}
\end{figure}
The challenge lies in determining the gradient, $\boldsymbol\nabla u|_{P_0}$ in equation \eqref{eqn:trialeq}, as it's typically unknown beforehand. To address this, FPM employs the Generalized Finite Difference (GFD) \cite{paper1} method that involves defining the support of the Point $P_0$ first. This is typically done by selecting a set of neighboring points surrounding the internal point $P_0$. The two common strategies for selecting neighbors are: (a) \textbf{Radius-based selection}, where all points $P \in \{P(\mathbf{x})| (\mathbf{x}-\mathbf{x_{0}})\leq r\}$, within a fixed radius $r$ are included, (b) \textbf{Nearest neighbor selection},  which considers points in subdomains sharing boundaries with the subdomain of $P_0$.  In this paper, we use the second approach, defining the support of $P_0$ as all points that share boundaries with the subdomain $E_0$, as shown in \Cref{fig:support}. These neighboring points are denoted as $P_1, P_2, P_3,....P_m$.  \\
\newline
\noindent The gradient $\boldsymbol{\nabla}u$ is determined by defining a weighted discrete $L^2$ norm $J$ \cite{paper1, liszka1980finite} as follows:
       \begin{align}
           J=\sum_{i=1}^m ((\mathbf{x_{i}}-\mathbf{x_{0}}).\boldsymbol\nabla{u}|_{P_0}-(u_{i}-u_{0}))^2 w_{i}
       \end{align}
where $\mathbf{x_{i}}$ denotes the coordinate vector of $P_{i}$, $u_{i}$ denotes the value of $u_{h}$ at $P_{i}$ and $w_{i}$ signifies the value of weight function of $P_0$ at $P_{i}~(i=1,2,...m)$. The weight function is assumed to remain constant in this context. The stationarity of $J$ leads to the following formula for the gradient $\boldsymbol\nabla u$ at ${P_0}$:
       \begin{align}\label{grad_u}
           \boldsymbol\nabla{u}|_{P_0}=\mathbf{B}\mathbf{u}_{E},
       \end{align}
 where the matrix $\mathbf{B}$ and the vector $\mathbf{u}_E$ are defined as
 \begin{subequations}
\begin{align}
\mathbf{B}&=(\mathbf{A}^{\rm T}\mathbf{A})^{-1} \mathbf{A}^{\rm T} \label{B_mat}\begin{bmatrix}\mathbf{I}_{1}&\mathbf{I}_{2}\end{bmatrix},\\
\mathbf{u}_{E}&=\left(\begin{bmatrix}u_{0}&u_{1}&...&u_{m}\end{bmatrix}\right)^{\rm T}.
\end{align}
\end{subequations}
The matrices $\mathbf{A}, \mathbf{I}_1$ and $\mathbf{I_2}$ are further defined as
\begin{subequations}
\begin{align}
\mathbf{A}&=\begin{bmatrix}
\mathbf{x}_{1}-\mathbf{x}_{0}\\
\mathbf{x}_{2}-\mathbf{x}_{0}\\
\cdots\\
\mathbf{x}_{m}-\mathbf{x}_{0}
\end{bmatrix},\\
\mathbf{I}_{1}&=\left(\begin{bmatrix}-1&-1&...&-1\end{bmatrix}_{1\times m}\right)^{\rm T},\\
\mathbf{I}_{2}&=\begin{bmatrix}
1&0&0&\cdots&0\\
            0&1&0&\cdots&0\\
            \cdots&\cdots&\cdots&\cdots&\cdots\\
            0&0&0&\cdots&1
            \end{bmatrix}_{m\times m}.
\end{align}
\end{subequations}
Finally, by substituting $\boldsymbol\nabla{u}|_{P_0}$ from equation \eqref{grad_u} into equation \eqref{eqn:trialeq}, the trial function $u_{h}$ becomes
       \begin{align}
        u_{h}(\mathbf{x})=\mathbf{N}\mathbf{u}_{E} , \qquad \mathbf{x}\in E_0
        \end{align}
where 
        \begin{equation}\label{shape_func}
            \mathbf{N}=(\mathbf{x}-\mathbf{x_{0}})\mathbf{B}+\begin{bmatrix}1&0&...&0\end{bmatrix}_{1\times(m+1)}.
        \end{equation}
The matrix $\mathbf{N}$ serves as the shape function of $u_h$ within the subdomain $E_0$, defined in terms of the internal point $P_0$ and its support points $P_1,P_2,P_3...,P_m$. Thus, neighbouring subdomains have their own shape function values at the common boundary, as no continuity requirement is imposed at the internal boundaries. Consequently, the shape functions are discontinuous at these boundaries. To demonstrate this, the shape function plots at a point in both 1-D and 2-D domains are shown in \Cref{fig:shpfunc1} and \Cref{fig:shpfunc2}, respectively.

\noindent Regarding the test function $v_h$, we stipulate that it adopts the identical shape function as $u_h$ within the weak Galerkin formulation. It is evident from the assembly of trial and test functions that they possess characteristics of being local, straightforward polynomials, reliant on point-based approaches, and exhibiting piecewise-continuous behavior throughout $\Omega$.

\section{Fragile Points formulation for governing equation}\label{numflux}
\subsection{Concept of point stiffness }
\begin{figure}[!htb]
    \centering
    \begin{minipage}{.45\textwidth}
    \includegraphics[scale=0.5]{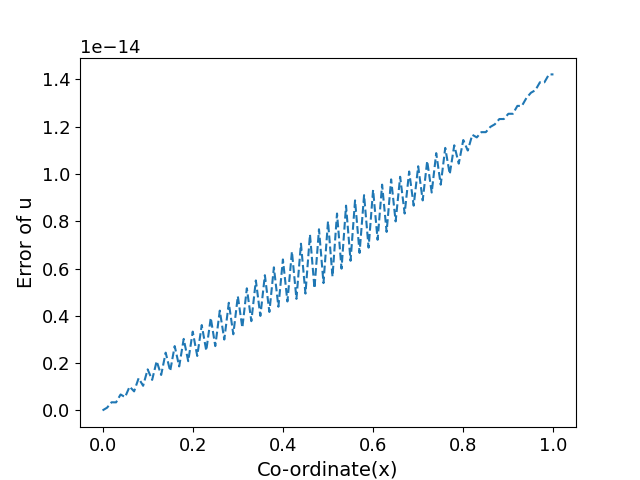}
    \captionsetup{justification=centering}
    \caption{ 101 uniformly distributed points}
    \label{fig:boundcond1}     
    \end{minipage}\hfill
\begin{minipage}{.45\textwidth}
    \includegraphics[scale=0.5]{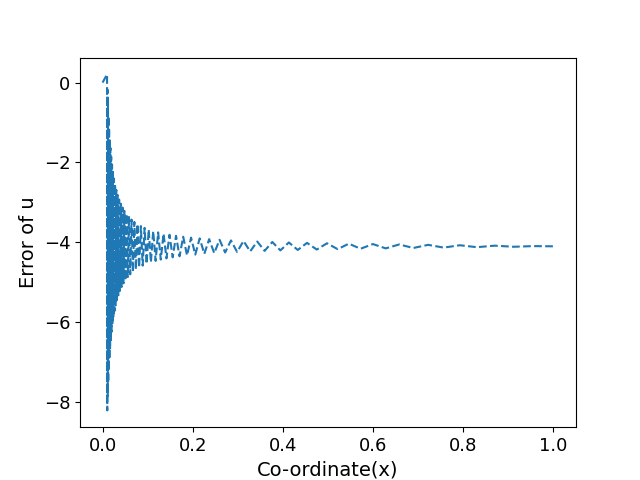}
    \captionsetup{justification=centering}
    \caption{ 101 randomly distributed points}
    \label{fig:boundcond2}
    \end{minipage}
    \end{figure}

\noindent In the governing \Cref{eqn:thermal_govern}, multiplying both sides by a test function $v$ and integrating by parts over the entire domain $\Omega$ yields the following weak form:
\begin{align}
    \sum\int_{E_i}\boldsymbol\nabla{v}^{T}\beta\boldsymbol\nabla u ~\rmd\Omega = \sum\int_{E_i}Q v ~\rmd\Omega + \int_{ \tilde{\Gamma}}v {\mathbf{n}}^{T}\beta \boldsymbol\nabla{u} ~\text{d}\Gamma.
    \label{eqn:weakform1}
\end{align}
where $\mathbf{n}$ denotes the outward unit normal vector on the boundary $\Gamma$. For a subdomain $E_0$, we substitute the trial and test functions $u_h$ and $v_h$ in terms of the shape function $\mathbf{N}$  and its gradient $\mathbf{B}$. Substituting these into the weak form \Cref{eqn:weakform1}, the point stiffness matrix $\mathbf{K_{P}}$ of the subdomain $E_0$ is derived as
\begin{align}
    \mathbf{K_{P}}=\int_{E_0}\mathbf{B}^{\rm T} \beta \mathbf{B} ~\rmd\Omega.
\end{align}
The use of linear interpolation for $u_h$ and $v_h$, combined with the computation of $\mathbf{B}$ exclusively at the internal point $P_0$ within $E_0$, ensures that the matrix $\mathbf{B}$ remains constant for this subdomain. As a result, when the diffusion coefficient $\beta$ is constant, the integration required for the calculation of the point stiffness matrix simplifies to
\begin{align}
    \mathbf{K_{P}}=\mathbf{B}^{\rm T}\beta\mathbf{B}A_{E_{0}},
\end{align}
where $A_{E_0}$ represents the area of the subdomain $E_0$. The global stiffness matrix for the entire domain $\Omega$ is then constructed by assembling the point stiffness matrices from all subdomains. Notably, unlike the finite element method (FEM) which uses element stiffness matrices, the fragile point method (FPM) employs point stiffness matrices. \answer{This is to note that although here, we use constant coefficients per cell, the formulation admits spatially varying, anisotropic tensors $\beta(\mathbf x)$ without requiring any modification.}\\
\newline
However, discontinuous trial and test functions in the Galerkin weak form lead to inconsistencies and inaccuracy in the solution failing the patch test. To demonstrate this, consider a one-dimensional steady-state heat conduction problem governed by the equation
\begin{subequations}
\begin{align}
    -\frac{\text{d}}{\text{d}x}\left(\beta\frac{\text{d}u}{\text{d}x}\right) &=1 ,\hspace{0.8cm}x\in (0,1),\\
    u_{x=0}&=0 ,\hspace{0.8cm} u_{x=1}=1.
\end{align}
\label{eq:1D}
\end{subequations}
The exact solution of \Cref{eq:1D} is given by $u=x$. \Cref{fig:boundcond1,fig:boundcond2} depict the error (difference between the exact and numerical solutions) obtained by solving \Cref{eq:1D} using FPM with 101 points distributed uniformly and randomly across the domain.
\noindent As evident from \Cref{fig:boundcond1}, the discontinuous trial and test functions lead to solution inconsistencies, particularly for randomly distributed points. To mitigate this issue, FPM incorporates numerical flux corrections inspired by the discontinuous Galerkin method, which will be discussed in the next subsection.

\subsection{Concept of numerical fluxes}
\label{sec:Concept of numerical fluxes}
As discussed in the previous section, the Fragile Points Method (FPM), which utilizes point-wise trial and test functions that are discontinuous at element boundaries, fails to deliver accurate results even with an increased number of points. To address this, numerical flux corrections are implemented to ensure consistency and enhance the method's accuracy. Commonly associated with the Discontinuous Galerkin method, numerical fluxes improve result accuracy and consistency. These fluxes facilitate the continuity of the solution across boundaries and contribute to stability. They typically incorporate jump and average operators for the primary and secondary variables. The average operators promote continuity at interfaces, while the jump operators help stabilize the solution.

To adapt the concept of numerical flux within the framework of the FPM, the governing differential equation (c.f. \Cref{eqn:thermal_govern}), is multiplied with the test function $v$ and using the Gauss divergence theorem, the weak form can be written as
\begin{align}
   \sum_{E\in \Omega}\int_{E}\boldsymbol\nabla{v}^{T}\beta\boldsymbol\nabla u ~\rmd\Omega= \sum_{E\in \Omega}\int_{E}Q v ~\rmd\Omega + \int_{ \tilde{\Gamma}}v \mathbf{n}^{T}\beta \boldsymbol\nabla{\hat{u}} ~\rmd\Gamma
    \label{eqn:weakform2}
\end{align}
where $\hat{u_{h}}$ is called the numerical flux, which is the approximation of $u_{h}$ on the boundary of subdomain $E$ for all $E \in \Omega$. Next, using the average($\{\}$) and the jump operator($\llbracket \rrbracket$), \Cref{eqn:weakform1} is rewritten as:
\begin{multline}
   \sum_{E\in \Omega}\int_{E}\boldsymbol\nabla{v}^{T}\beta\boldsymbol\nabla u ~\rmd\Omega= \sum_{E\in \Omega}\int_{E}Q v ~\rmd\Omega+\sum_{e\in \Gamma_{h}}\int_{e}\left(\{v\}\llbracket \beta\boldsymbol\nabla{\hat{u}} \rrbracket+\llbracket v \rrbracket\cdot \{\beta\boldsymbol\nabla{\hat{u}}\}\right)d\Gamma\\ + \sum_{e\in \Gamma_{D}\bigcup\Gamma_{N}}\int_{e}\llbracket v \rrbracket \cdot\{\beta\boldsymbol\nabla{\hat{u}}\}d\Gamma
    \label{eqn:weakform_with jump}
\end{multline}
where $\tilde{\Gamma} = \Gamma_h \cup \Gamma_D \cup \Gamma_N$ is the set of all boundaries, $\Gamma_h$ is the set of all internal boundaries and Dirichlet and Neumann conditions are specified on $\Gamma_{D}$, $\Gamma_{N}$, respectively, representing the external boundaries. The average and jump operators for a scalar quantity $\phi$ and a vector quantity 
$\boldsymbol{\Phi}$ are defined as follows:
\begin{linenomath}
\begin{subequations}
\begin{align}
& \forall e\in \Gamma_{h} \left\{         \begin{aligned}
             \{\phi\}&=\frac{1}{2}(\phi_{1}+\phi_{2}),\quad \llbracket \phi\rrbracket=\phi_{1}\mathbf{n}_{1}+\phi_{2}\mathbf{n}_{2}\\ 
             \{\boldsymbol{\Phi}\}&=\frac{1}{2}(\boldsymbol{\Phi}_{1}+\boldsymbol{\Phi}_{2}),\quad \llbracket\boldsymbol{\Phi}\rrbracket=\boldsymbol{\Phi}_{1}\cdot\mathbf{n}_{1}+\boldsymbol{\Phi}_{2}\cdot\mathbf{n}_{2} \end{aligned}
             \right. \\
& \forall e \in [\Gamma_{D}, \Gamma_{N}] \quad 
         \{\phi\}=\phi,\quad \llbracket \phi\rrbracket=\phi\mathbf{n},\quad \{\boldsymbol{\Phi}\}=\boldsymbol{\Phi},\quad \llbracket\boldsymbol{\Phi}\rrbracket=\boldsymbol{\Phi}\cdot\mathbf{n}     
             \end{align}
\end{subequations}
\end{linenomath}

\subsection{Interior Penalty (IP) Numerical Flux and Primal form}
There are different types of Numerical Fluxes that have evolved and been studied by many authors (see the review by \cite{paper2}). Proper choice of numerical flux can lead to a stable, accurate, and consistent solution and also affect the computational time of a particular method. In this paper, we use Interior Penalty (IP) as the Numerical Flux as shown in \Cref{T:table1}.
\begin{table}
\caption{Interior Penalty Numerical Flux}
\centering
\begin{tabular}{cccc}
\hline
  Keyword &$\Gamma_{h}$ &$\Gamma_{D}$&$\Gamma_{N}$ \\
\hline
 $\hat{u}$&$\{u\}$& $g_{D}$&$u$\\
$\boldsymbol\nabla{\hat{u}}$ & $\{\boldsymbol\nabla u\}-\frac{\eta}{h_{e}}\llbracket u\rrbracket$ &$\boldsymbol\nabla u-\frac{\eta}{h_{e}}(u-g_{D})\mathbf{n}$& $\boldsymbol{g}$ \\
\hline
\end{tabular}
\label{T:table1}
\end{table}
In \Cref{T:table1}, $h_{e}$ is a boundary dependent parameter with the dimensions of length and can be defined as length of the subdomain in case  of 1D and length of the boundary in 2D, $\eta$ is the penalty parameter (positive number) independent of the edge size. Upon substituting the IP numerical flux into \Cref{eqn:weakform_with jump} we get the primal form of the FPM:
\begin{multline}
   \sum_{E\in \Omega}\int_{E}\boldsymbol\nabla{v}^{T}\beta\boldsymbol\nabla u ~\rmd\Omega-\sum_{e\in \Gamma_{h}\cup\Gamma_{D}}\int_{e}\left(\{\beta\boldsymbol\nabla{u}\}\cdot\llbracket v \rrbracket+\{\boldsymbol\nabla{v}\}\cdot\llbracket\beta u \rrbracket \right)\text{d}\Gamma+\sum_{e\in \Gamma_{h}\cup\Gamma_{D}}\frac{\eta}{h_{e}}\int_{e}\llbracket \beta u\rrbracket\cdot\llbracket v\rrbracket ~\text{d}\Gamma\\=\sum_{E\in \Omega}\int_{E}Q v ~\rmd\Omega+\sum_{e\in\Gamma_{D}}\int_{e}\beta\left(\frac{\eta}{h_{e}}v-\boldsymbol\nabla{v}\cdot \mathbf{n}\right)g_{D}~\text{d}\Gamma+\sum_{e\in\Gamma_{N}}\int_{e}v\beta g_N~\text{d}\Gamma
    \label{eq:primalform}
\end{multline}

\subsection{Numerical Implementation}
The primal form of FPM \Cref{eq:primalform} can be written in the matrix form:
\begin{equation}\label{final_eq}
\mathbf{K}\mathbf{u} = \mathbf{f}
\end{equation}
where $\mathbf{K}$ is the global thermal conductivity matrix, $\mathbf{u}$ is the unknown vector with nodal DOFs and $\mathbf{f}$ is the load vector. Upon substituting $\mathbf{N}$ for $u$ and $v$, $\mathbf{B}$  for $\boldsymbol\nabla u$ and $\boldsymbol\nabla v$ in \Cref{eq:primalform}, the point stiffness matrix $\mathbf{K_{P}}$, point stiffness matrix at internal boundary $\mathbf{K_{IP}}$, Dirichlet boundary stiffness matrix $\mathbf{K_{D}}$ and load matrix $\mathbf{f}$ can obtained as :
\begin{subequations}
    \begin{align}
        \mathbf{K_{P}} &=\int_{E}\mathbf{B}^{\rm T} \beta \mathbf{B}~\rmd\Omega, \qquad\text{where}\quad E \in \Omega \label{KP}\\ 
         \begin{split}\hspace{-1cm}
             \mathbf{K_{IP}}=-\frac{1}{2}\int_{e}(\beta_1\mathbf{B_{1}}^{\rm T} \mathbf{n}^{\rm T}\mathbf{N_{1}}+\beta_1\mathbf{N_{1}}^{\rm T} \mathbf{n}  \mathbf{B_{1}})~\rmd\Gamma+\frac{\eta}{h_{e}}\int_{e}\beta_1\mathbf{N_{1}}^{\rm T}\mathbf{N_{1}}~\rmd\Gamma\\-\frac{1}{2}\int_{e}(\beta_2\mathbf{B_{2}}^{\rm T} \mathbf{n}^{\rm T}\mathbf{N_{2}}+\beta_2\mathbf{N_{2}}^{\rm T}\mathbf{n} \mathbf{B_{2}})~\rmd\Gamma+\frac{\eta}{h_{e}}\int_{e}\beta_2\mathbf{N_{2}}^{\rm T}\mathbf{N_{2}}~\rmd\Gamma\\-\frac{1}{2}\int_{e}(\beta_2\mathbf{B_{2}}^{\rm T} \mathbf{n}^{\rm T}\mathbf{N_{1}}+\beta_2\mathbf{N_{2}}^{\rm T}\mathbf{n}  \mathbf{B_{1}})~\rmd\Gamma+\frac{\eta}{h_{e}}\int_{e}\beta_1\mathbf{N_{1}}^{\rm T}\mathbf{N_{2}}~\rmd\Gamma\\-\frac{1}{2}\int_{e}(\beta_1\mathbf{B_{1}}^{\rm T}  \mathbf{n}^{\rm T}\mathbf{N_{2}}+\beta_1\mathbf{N_{1}}^{\rm T}\mathbf{n}  \mathbf{B_{2}})~\rmd\Gamma+\frac{\eta}{h_{e}}\int_{e}\beta_2\mathbf{N_{2}}^{\rm T}\mathbf{N_{1}}~\rmd\Gamma,
             \end{split}\label{K_IP}
             \\
             &\text{where}\quad e =\partial E_1\cap \partial E_2\nonumber \\
            \mathbf{K_{D}}&=\frac{1}{2}\int_{e}(\beta\mathbf{B}^{\rm T}  \mathbf{n}^{\rm T}\mathbf{N}+\beta\mathbf{N}^{\rm T}\mathbf{n}  \mathbf{B})~\rmd\Gamma+\frac{\eta}{h_{e}}\int_{e}\beta\mathbf{N}^{\rm T}\mathbf{N}~\rmd\Gamma,\quad \text{where}\quad e\in\Gamma_D\label{KD}\\
             \mathbf{f}&=\int_{E}Q\mathbf{N}~\rmd\Omega+\int_{\Gamma_{D}}\beta\Big(\frac{\eta}{h_{e}}\mathbf{N}- \mathbf{B}\Big)g_{D}~\rmd\Gamma+\int_{\Gamma_{N}}\mathbf{N}\beta g_{N}~\rmd\Gamma\label{f}
\end{align}
\end{subequations}

The total stiffness matrix is then
given by:
\begin{equation}
    \mathbf{K}=\mathbf{K_{P}}+\mathbf{K_{IP}}+\mathbf{K_{D}}.
\end{equation}
Since the trial function used in FPM is discontinuous, interface problems can be modeled explicitly within this method by means of the following simple modifications. When an interior interface is there, the support of the pair of neighbouring points should be modified to exclude each other from their support, which guarantees that the trial functions in this pair of neighbouring subdomains become independent. In this way, the interaction between the neighbouring subdomains is completely separated when the interface is formed. A pseudo-code for the implementation of the fragile points method is given in \Cref{algorithm}.

\begin{algorithm}[htpb]
\DontPrintSemicolon
\caption{Algorithm to solve Interface problem}
\label{algorithm}
\KwInput{Coordinates of the fragile points}
\KwOutput{Stiffness matrix}
\For{i = 1:{\rm number of fragile points}}
{
Compute shape functions $\mathbf{N}$ from Eq.~\eqref{shape_func} \;
Compute $\mathbf{B}$ using Eq. \eqref{B_mat} \;
Compute the point stiffness matrix $\mathbf{K_{P}}$ using Eq. \eqref{KP} and assemble to the global stiffness matrix, $\mathbf{K}$ \;
\If{the interface involves a primary variable jump}
{
During assemble of point stiffness to global stiffness matrix, don't connect the two side interface elements
}
\Else
{
Compute the Global stiffness matrix as we do for the finite element method
}
}
\For{e = 1:{\rm number of internal interface edges}}
{
Identify subdomains $E_1$ and $E_2$ sharing the $e$-th interface edge \;
Compute shape functions $\mathbf{N_1},~\mathbf{N_2}$ on $E_1$ and $E_2$ \;
Compute gradients $\mathbf{B_1},~\mathbf{B_2}$ \;
Compute the interface stiffness matrix $\mathbf{K_{IP}}$ using Eq.~\eqref{K_IP} \;
}
Compute the Dirichlet boundary stiffness matrix $\mathbf{K_D}$ using Eq.~\eqref{KD} \;
Assemble all $\mathbf{K_{P}}$, $\mathbf{K_{IP}}$ and $\mathbf{K_D}$ to find the global stiffness matrix $\mathbf{K}$ \;
Compute the load vector $\mathbf{f}$ using Eq. \eqref{f} \;
Implement boundary conditions either strongly or weakly \;
Solve the assembled Eq.~\eqref{final_eq} to find the primary variable values \;
\end{algorithm}

 \section{Numerical examples}
\label{numexap}
This section evaluates the accuracy and convergence properties of the proposed framework by applying it to various benchmark problems. 
We consider example problems with interfaces of various shapes, including straight, circular and complex star-shaped interfaces. In order to estimate the errors and study the convergence properties of the numerical results, the following relative error ($r$) with the displacement $L^2$ norm is employed: 
\begin{equation}\label{rel_err_l2nrm}
    \begin{aligned}
    r =\frac{||\mathbf{u}^h - \mathbf{u}^{exact}||_{L^2}}{||\mathbf{u}^{exact}||_{L^2}}
\end{aligned}
\end{equation}
where,
\begin{align*}
    ||\mathbf{u}||_{L^2}=\sqrt{\int_{\Omega}\mathbf{u^{\rm T}~u}~\rmd\Omega}
\end{align*}
and error (\%) as
\begin{align*}
    error (\%) =\left|\frac{\mathbf{u}^h - \mathbf{u}^{exact}}{\mathbf{u}^{exact}}\right|\times100.
\end{align*}
For the Fragile Point Method (FPM) simulations, linear trial and test functions are utilized. The paper also showcases a comparison between the FPM simulations that use numerical fluxes and the analytical solutions for each example in the subsequent subsections.

\subsection{Patch test}

\begin{figure}[!htp]
\centering
\begin{subfigure}[b]{0.5\textwidth}
\centering
\includegraphics[width=1.1\textwidth]{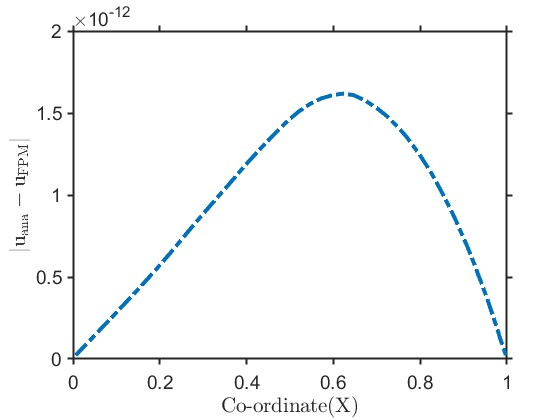}
\caption{}
\label{sf:uniform(with NF)}
\end{subfigure}\hfill
\begin{subfigure}[b]{0.5\textwidth}
\centering
\includegraphics[width=1.1\textwidth]{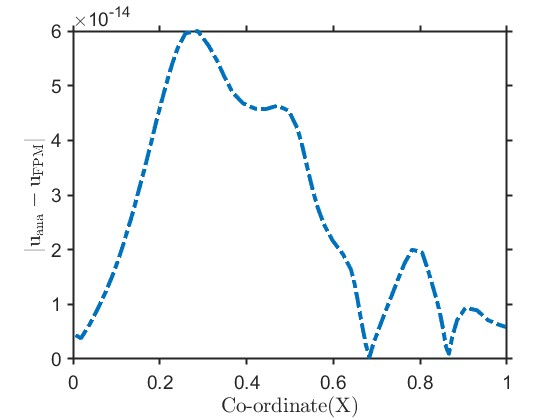}
\caption{}
\label{sf:random(with NF)}
\end{subfigure}
\caption{(a) 1D-101 uniformly distributed points (with NF) (b) 1D-101 randomly distributed points (with NF) }
\label{fig:boundcond patch 1D_thermal}
\end{figure}

To examine the consistency of the Fragile Points Method (FPM), we conducted a 1D patch test to solve the one-dimensional steady-state heat conduction equation governed by the following equation:
\begin{equation}
    \frac{\text{d}^2u}{\text{d}x^2}=0, \qquad x\in[0,1].
\end{equation}
 The analytical solution for this patch test is $u(x)=x$, representing a linear temperature gradient across the domain. Our one-dimensional domain extends from $x=0$ to $x=1$. The boundary conditions are set such that $u(0)=0$ and $u(1)=1$, perfectly aligning with the chosen analytical solution.
\newline
\indent During the FPM simulation, linear trial and test functions are employed, utilizing a penalty coefficient $\eta=10$ for this 1D scenario. Since the exact solutions are linear functions, to pass the patch tests, the numerical solutions must equal the exact solutions up to machine precision. The error in $u$ for the 1D case, for both uniformly distributed and randomly distributed points, is illustrated in \Cref{fig:boundcond patch 1D_thermal}. This comparison infers that, due to the implementation of numerical flux during the calculations, our proposed methodology is consistent in both scenarios, achieving solutions that are accurate to machine precision.
\subsection{Straight interface}
\begin{figure}[!htb]
 \subfloat[]{\includegraphics[scale=0.35]{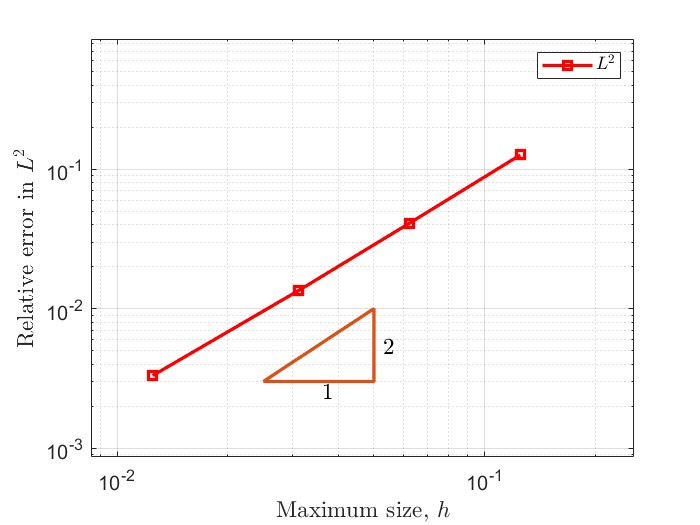}\label{fig:straight_interface_wo_CR}}
\subfloat[]{\includegraphics[scale=0.35]{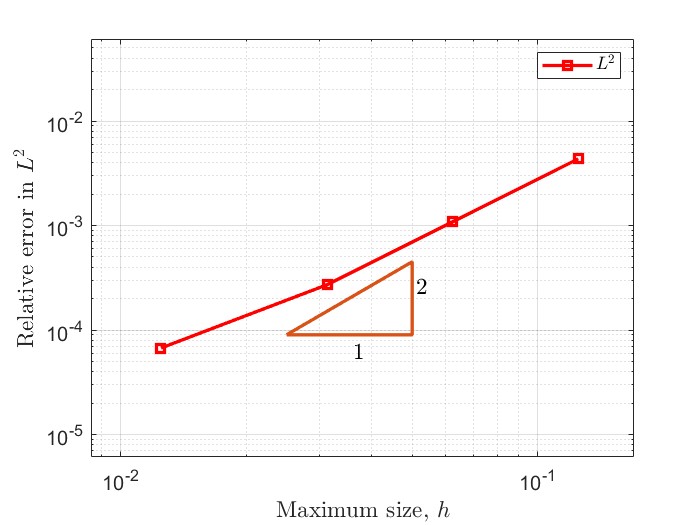}\label{fig:straight_interface_w_CR}}
\caption{Convergence of the relative error in the $L_{2}$ norm for a domain with the straight interface: (a) with zero jump condition at the interface, i.e., $\llbracket u \rrbracket=0$ and (b) with a non-zero jump condition at the interface, i.e., $\llbracket u \rrbracket=\delta$.}
\end{figure}

Consider a one-dimensional diffusion problem within a two-dimensional domain denoted as $\Omega = [0,1]\times[0,1]$. Here, the domain $\Omega$ is divided into two distinct subdomains by a straight interface located at $x=0.5$. The subdomain $\Omega_1$ is characterized by thermal conductivity $\beta_1$ and $\Omega_2$ is characterized by the thermal conductivity $\beta_2$. 
As shown in \Cref{fig:straight_interface_domain}(a), zero Neumann boundary conditions are applied at the top $ (y = 1)$ and at the bottom $(y = 0)$ boundaries, signifying that these surfaces are perfectly insulated. In contrast, Dirichlet boundary conditions are employed at the left $(x = 0)$ and the right $(x = 1)$ boundaries, where the temperature is specified from the analytical solution given by \cite{dsouza2021robust}. At the interface, the jump conditions are specified by
\begin{subequations}
    \begin{align}
    \llbracket u \rrbracket&=\delta, \qquad
    \llbracket \beta\boldsymbol\nabla u \rrbracket=0,
\end{align}
\end{subequations}
where $\delta$ represents the prescribed jump in temperature. Two scenarios are examined: a zero jump condition ($\delta= 0$), indicating continuous temperature across the interface, and a non-zero jump condition ($\delta = 2$), indicating a discontinuity in temperature. The thermal conductivities are selected such that the ratio $\dfrac{\beta_{1}}{\beta_{2}}= 6$. Additionally, a constant heat source term, $Q = 1$, is assumed throughout the domain. The analytical solutions for these boundary conditions are as follows~\cite{dsouza2021robust}:
\begin{figure}[!htb]
\subfloat[]{\includegraphics[scale=0.73]{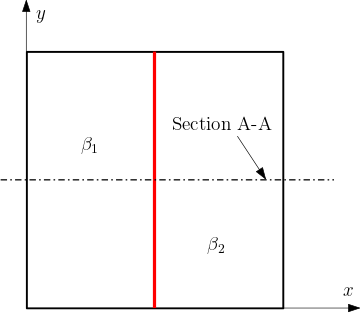}}~~
\subfloat[]{\includegraphics[scale=0.48]{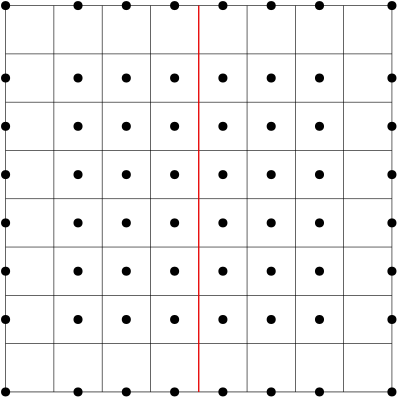}\label{fig:strmesh}}
\caption{(a) Square domain with straight interface geometry and its boundary conditions and (b) representative mesh used for FPM}
\label{fig:straight_interface_domain}
\end{figure}
\begin{subequations}
    \begin{align}
    u_{1}(x)&=\dfrac{(3\beta_{1}+\beta_{2})x}{4\beta_{1}^2+4\beta_{1}\beta_{2}}-\dfrac{x^2}{2\beta_{1}} \quad\quad\quad\quad\quad\quad\quad \quad\quad \forall x<0.5\\
    u_{2}(x)&=\dfrac{\beta_{2}-\beta_{1}+(3\beta_{1}+\beta_{2})x}{4\beta_{1}^2+4\beta_{1}\beta_{2}}-\dfrac{x^2}{2\beta_{2}}+\delta \quad\quad\quad \forall x\geq0.5
\end{align}
\end{subequations}
\begin{figure}[!htb]
 \subfloat[]{\includegraphics[scale=0.55]{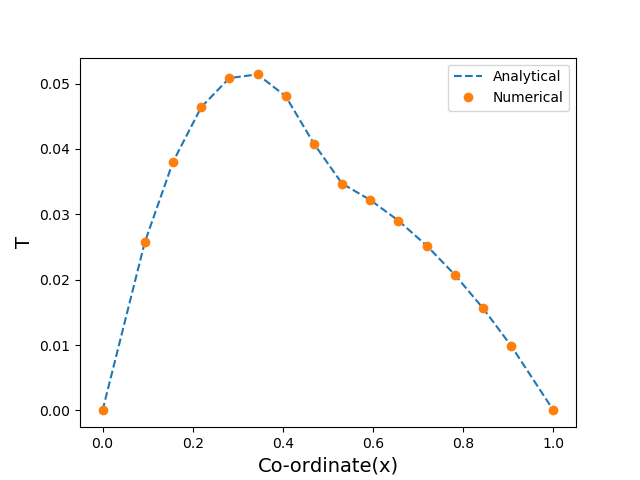}\label{fig:straight_interface_results_wo_jump}}
\subfloat[]{\includegraphics[scale=0.55]{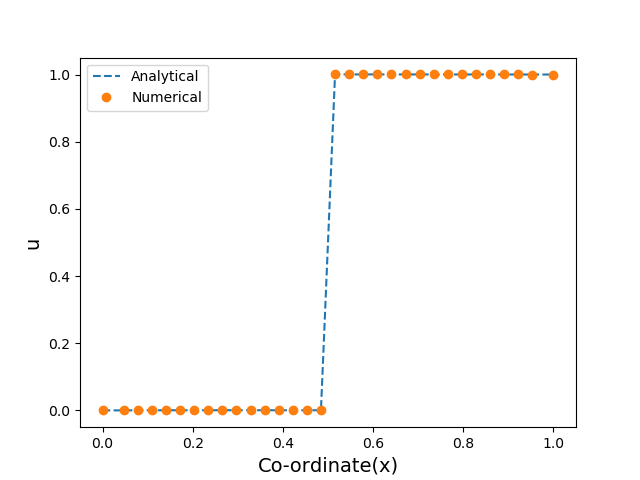}\label{fig:straight_interface_results_w_jump}}
\caption{(a) Numerical solution of the temperature over the domain from the proposed framework for without jump, (b) Numerical solution of the temperature over the domain from the proposed framework for with jump at interface}
\label{fig:straight_interface_results}
\end{figure}
For numerical simulations, the domain is refined with a total of 64 points, as shown in \Cref{fig:strmesh}. A penalty coefficient of $\eta=$ 1 was used. The resulting temperature profiles along section A-A within the domain are depicted in \Cref{fig:straight_interface_results} wherein \Cref{fig:straight_interface_results_wo_jump} showcases the temperature distribution for the zero jump condition scenario, while \Cref{fig:straight_interface_results_w_jump} presents the results for the non-zero jump condition scenario. These results show a good agreement with the analytical solutions. Furthermore, to assess the accuracy of the numerical solutions as the mesh is refined, we analyze the convergence of the relative error in the $L_2$ norm using \Cref{rel_err_l2nrm}. The convergence behavior for both the zero jump and nonzero jump scenarios is illustrated in \Cref{fig:straight_interface_wo_CR,fig:straight_interface_w_CR}, respectively. It is seen that the results from the FPM is accurate and converges at a rate of 2.

\subsection{Circular interface}
\begin{figure}[!htb]
\centering
\includegraphics[scale=0.7]{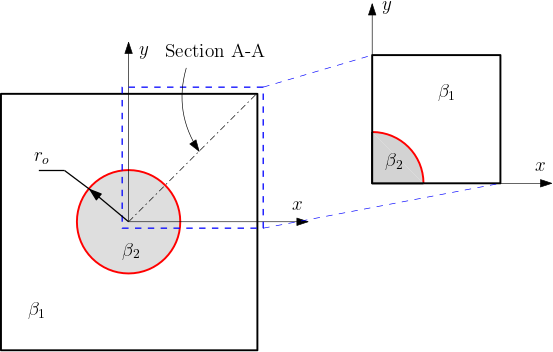}
\caption{Square domain with circular interface geometry and its boundary conditions and the simplified quarter model after considering symmetry}
\label{fig:circular_interface_domain}
\end{figure}

\begin{figure}[!htb]
 \subfloat[]{\includegraphics[scale=0.35]{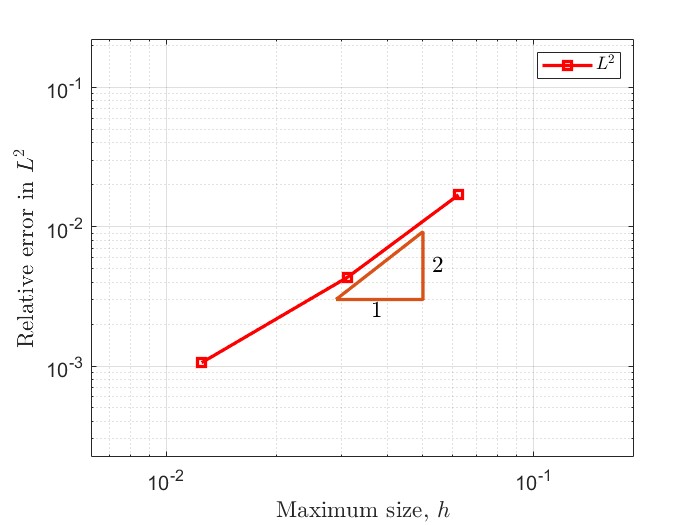}\label{fig:circ_interface_wo_CR}}
\subfloat[]{\includegraphics[scale=0.35]{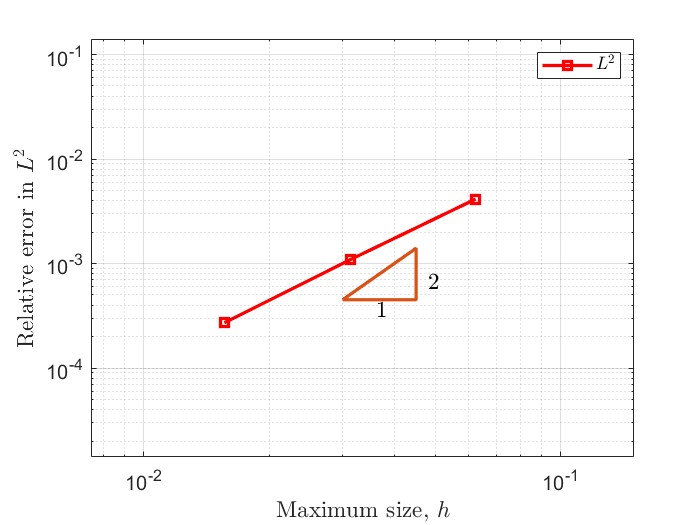}\label{fig:circ_interface_w_CR}}
\caption{Relative error in the $L_{2}$ norm for the circular interface: (a) without interface jump, $\llbracket k\boldsymbol\nabla u \rrbracket=0$ and (b) with $\llbracket k\boldsymbol\nabla u \rrbracket=\delta$.}
\end{figure}

This subsection explores a less trivial two-dimensional heat conduction problem within a square domain denoted by $\Omega = [-1,1]\times[-1,1]$. Heat conduction is governed by \Cref{eqn:thermal_govern}. The domain $\Omega$ encompasses a circular inclusion of radius $r_0 =$ 0.5 centered at its geometric centre, as illustrated in \Cref{fig:circular_interface_domain}. The thermal conductivities within the subdomains $\Omega_1$ and $\Omega_2$ are assigned the values of $\beta_1 = 1 $ and $\beta_2 = 1000$, respectively. In addition, a constant heat source term of $Q = -4$ is introduced throughout the domain. Due to the inherent symmetry in both the geometry and the boundary conditions, it suffices to solve for the temperature distribution in only one quadrant of the domain, see \Cref{fig:circular_interface_domain} for the geometry. The jump conditions at the interface are given by:
\begin{subequations}
    \begin{align}
    \llbracket u \rrbracket&=0 \qquad \text{and} \qquad
    \llbracket \beta\boldsymbol\nabla u \rrbracket=\delta,
\end{align}
\end{subequations}
where $\delta$ is the specified jump.
Similar to the previous example, two different cases are considered: a zero jump condition $(\delta=0)$, in the heat flux at the interface and a nonzero jump condition $(\delta=2)$,  in the heat flux at the interface. For the given set of boundary conditions, the analytical solution for the zero jump case is given by \cite{dsouza2021robust}:
\begin{subequations}
    \begin{align}
    u_{1}(x,y)&=\dfrac{r^2}{\beta_{1}} \quad\quad\quad\quad\quad\quad\quad \quad \forall r\leq0.5\\
    u_{2}(x,y)&=\dfrac{r^2}{\beta_{2}}-\dfrac{r_{0}^2}{\beta_{2}}+\dfrac{r_{0}^2}{\beta_{1}} \quad\quad\quad \forall r>0.5
\end{align}
\end{subequations}
and for the non-zero jump condition, the analytical solution is given by \cite{dsouza2021robust}:
\begin{subequations}
    \begin{align}
    u_{1}(x,y)&=1 \quad\quad\quad\quad\quad\quad \forall r\leq0.5\\
    u_{2}(x,y)&=1+ ln 2r \quad\quad\quad \forall r>0.5
\end{align}
\end{subequations}
\begin{figure}[!htb]
 \subfloat[]{\includegraphics[scale=0.55]{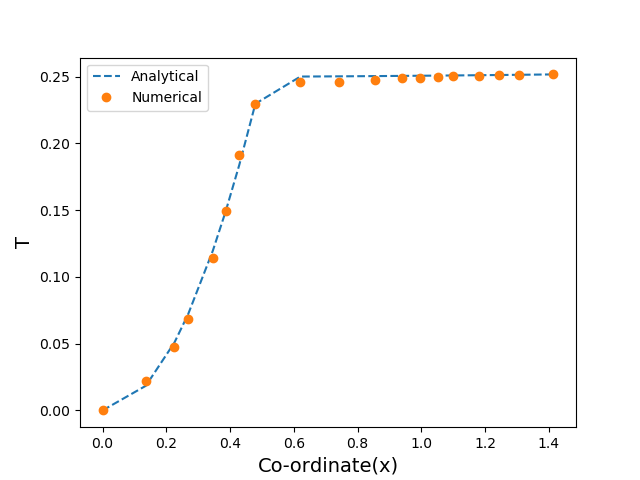}\label{fig:circ_interface_results_wo_jump}}
\subfloat[]{\includegraphics[scale=0.55]{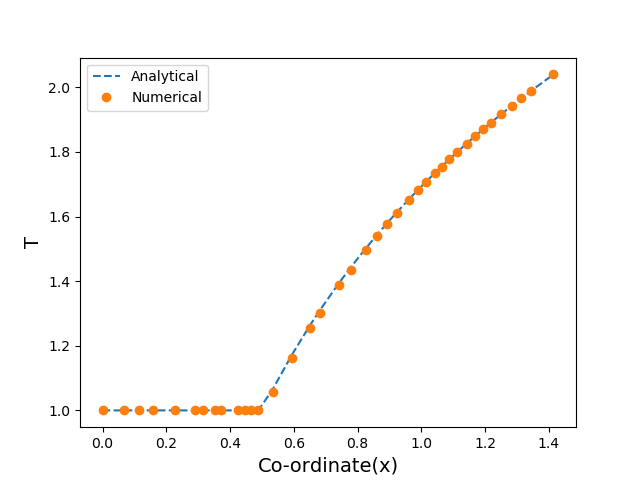}\label{fig:circ_interface_results_w_jump}}
\caption{(a) Numerical solution of the temperature over the domain from the proposed framework for without jump at the circular interface, (b) Numerical solution of the temperature over the domain from the proposed framework for with jump in the derivative at the circular interface}
\label{fig:circular_interface_result}
\end{figure}
where $r=\sqrt{x^2+y^2}$. For numerical simulation, a penalty coefficient of $\eta = 10$ was employed. The resulting temperature profiles along section A-A within the domain are presented in \Cref{fig:circular_interface_result} wherein \Cref{fig:circ_interface_results_wo_jump} showcases the temperature distribution for the zero jump condition, while \Cref{fig:circ_interface_results_w_jump} presents the results for the non-zero jump condition. As observed, the numerical solutions obtained from the proposed framework exhibit good agreement with the corresponding analytical solutions. The convergence of the relative error in the $L_{2}$ norm with mesh reﬁnement for both the cases is shown in \Cref{fig:circ_interface_wo_CR,fig:circ_interface_w_CR}.  It is seen that the results from the FPM is accurate and converges at a rate of 2.
\subsection{Star shaped interface}
\begin{figure}[htpb]
\subfloat[]{\includegraphics[scale=0.65]{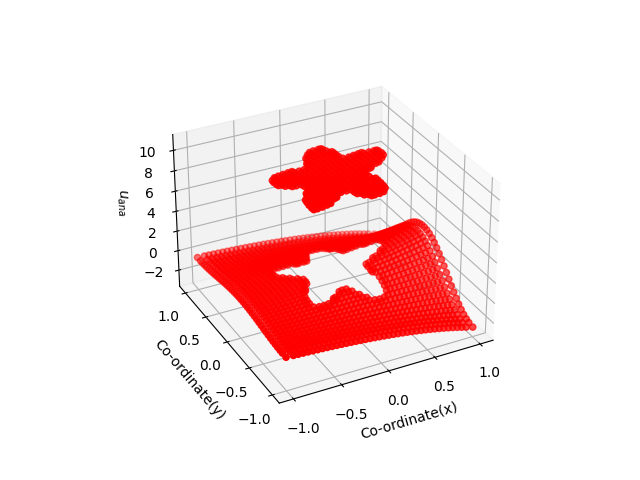}}
\subfloat[]{\includegraphics[scale=0.65]{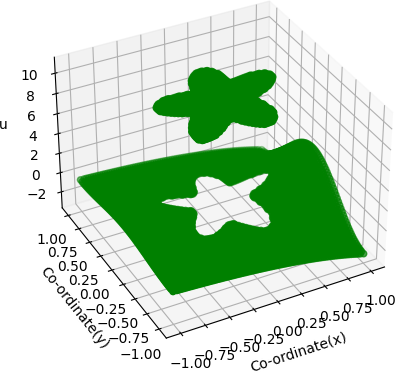}} 
\caption{(a)Numerical solution of the temperature over the domain from the proposed framework for with jump in the primary variable at the star-shaped interface (b) Analytical solution of the temperature over the domain from the proposed framework for with jump in the primary variable at the star-shaped interface}
\label{fig:star_interface}
\end{figure}

We further examine the convergence of the proposed method by considering a star shaped inclusion within a two-dimensional square domain $\Omega = [-1,1]\times[-1,1]$ as shown in \Cref{fig:star_interface}. The interface is generated by using the following equation in polar form:
\begin{align*}
    r = 0.45 + 0.1(\cos(5\theta) + \sin(5\theta)),\quad \text{with} \quad\theta
    \in [0, 2\pi).
\end{align*}
The thermal conductivity in the interior of the interface is assigned the value of $\beta_{1} = 100$ and on the exterior region, the value assigned is $\beta_{2} = 1$. The constant heat source term is $Q = 1$ throughout the domain and Dirichlet boundary conditions are imposed at the $x=0$, $x=1$, $y=0$ and $y=1$. The jump conditions at the interface are given by
\begin{equation}
    \llbracket u \rrbracket =\delta, \qquad
    \llbracket \beta\boldsymbol\nabla u \rrbracket = 0,
\end{equation}
where $\delta$ is the specified jump. For $\delta > 0$, non-zero jump condition in the temperature field at the interface the analytical solution is given by \cite{mu2016new}:
    \begin{align}
    u(\xx)=\begin{dcases}
        x-y^{2}+10, \quad\quad & \forall \quad \xx \in \Omega_{1},\\
        e^{x}\cos(\pi y) \quad\quad & \forall \quad \xx\in \Omega_{2}.
    \end{dcases}
\end{align}
\begin{figure}[htpb]
\centering
\includegraphics[width=0.55\linewidth]{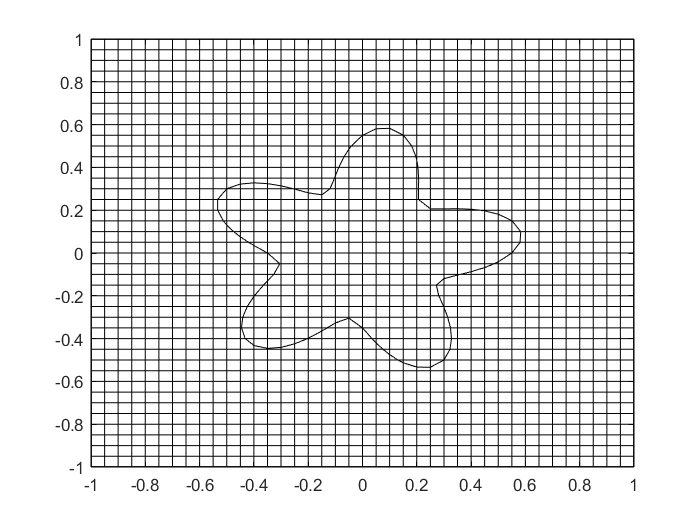}
\caption{A square domain mesh with the star-shaped interface}
\label{fig:figg1}
\end{figure}

\begin{figure}[htpb]
\centering
\includegraphics[width=0.75\linewidth]{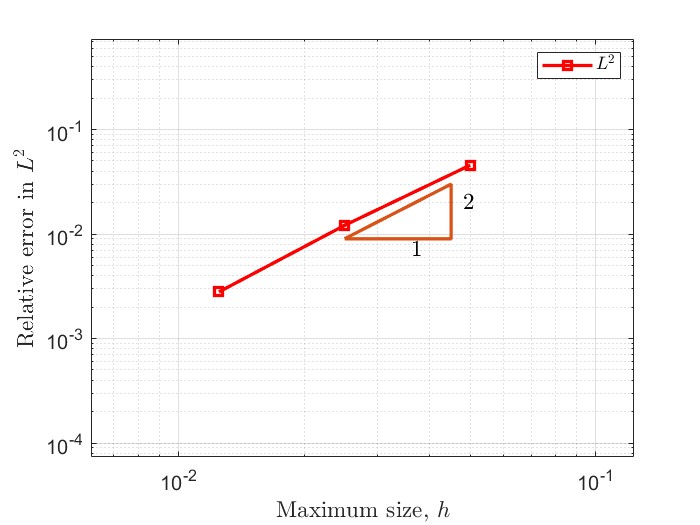}
\caption{Relative error in $L^2$ norm for the star interface}
\label{fig:starInt_w_jump}
\end{figure}
For numerical simulation, a penalty coefficient of $\eta =$ 10 was employed. The numerical and analytical temperature distributions over the domain for the non-zero jump condition in the primary variable at the star-shaped interface are illustrated in \Cref{fig:star_interface}, and a close agreement between the two can be observed. The convergence of the relative error in the $L_{2}$ norm with mesh reﬁnement for both the cases is shown in \Cref{fig:starInt_w_jump}.  It is seen that the results from the FPM is accurate and converges at a rate of 2, as expected.

\section{Conclusions}
This study demonstrates the effectiveness of the Fragile Points Method (FPM) as a robust and accurate meshless framework for solving two-dimensional interface problems in heterogeneous continua. By leveraging discontinuous shape functions constructed through a generalized finite difference approach, the method naturally enforces jump conditions on both primary and secondary variables without requiring additional constraints. The implicit representation of interfaces via level sets enhances the flexibility of the framework, which is validated across a range of benchmark problems with complex interface geometries.\\

\answer{Numerical results confirm that the proposed approach achieves high accuracy while maintaining optimal convergence rates in respective norms. The same flux-based, meshless framework could be extended to functionally graded porous media with TPMS architectures by representing the surface implicitly and embedding spatially varying transport tensors in the cell integrals. This positions FPM as a practical tool for graded heat or flow problems \citep{Tran2024APM,Tran2025METS}. Overall, this framework offers a promising tool for modeling interface phenomena, with potential extensions to moving boundary problems and three-dimensional applications planned for future research.}

\section*{Conflict of Interest}
The authors declare that they have no conflicts of interest regarding the publication of this manuscript.

\bibliographystyle{unsrt}
\bibliography{myRef}
\end{document}